\documentclass[twoside,12pt]{article}
\newcommand{\chapter}{\section}

\begin{document}

\newtheorem{Thm}{Theorem}
\newtheorem{Ax}{Axiom}
\newtheorem{Prop}{Proposition}
\newtheorem{Cor}[Prop]{Corollary}
\newtheorem{Main}{}
\renewcommand{\theMain}{}
\newtheorem{Lem}[Prop]{Lemma}
\newtheorem{Fact}{Fact}
\renewcommand{\theFact}{}

\newtheorem{Def}{Definition}
\newtheorem{rmk}{Remark}
\newenvironment{Rmk}{\begin{rmk}\em}{\end{rmk}}
\newtheorem{exm}{Example}
\newenvironment{Exm}{\begin{exm}\em}{\end{exm}}

\newcommand{\qed}{\par {\EM QED} }
\newtheorem{prf}{Proof}
\renewcommand{\theprf}{}
\newenvironment{Prf}{\begin{prf}\em}{\qed\end{prf}}
\newtheorem{prff}{}
\renewcommand{\theprff}{}
\newenvironment{Prff}{\begin{prff}\em}{\qed\end{prff}}



\newcommand{\YES}[1]{#1}
\newcommand{\NOT}[1]{}

\newcommand{\cA}{{\cal A}}
\newcommand{\cB}{{\cal B}}
\newcommand{\cC}{{\cal C}}
\newcommand{\cD}{{\cal D}}
\newcommand{\cE}{{\cal E}}
\newcommand{\cF}{{\cal F}}
\newcommand{\cG}{{\cal G}}
\newcommand{\cH}{{\cal H}}
\newcommand{\cI}{{\cal I}}
\newcommand{\cJ}{{\cal J}}
\newcommand{\cK}{{\cal K}}
\newcommand{\cL}{{\cal L}}
\newcommand{\cM}{{\cal M}}
\newcommand{\cN}{{\cal N}}
\newcommand{\cO}{{\cal O}}
\newcommand{\cP}{{\cal P}}
\newcommand{\cQ}{{\cal Q}}
\newcommand{\cR}{{\cal R}}
\newcommand{\cS}{{\cal S}}
\newcommand{\cT}{{\cal T}}
\newcommand{\cU}{{\cal U}}
\newcommand{\cV}{{\cal V}}
\newcommand{\cW}{{\cal W}}
\newcommand{\cX}{{\cal X}}
\newcommand{\cY}{{\cal Y}}
\newcommand{\cZ}{{\cal Z}}

\newcommand{\bbb}[1]{{\mbox{\bf #1}}}

\newcommand{\bN}{\bbb{N}}
\newcommand{\bZ}{\bbb{Z}}
\newcommand{\bR}{\bbb{R}}
\newcommand{\bC}{\bbb{C}}
\newcommand{\bQ}{\bbb{Q}}
\newcommand{\bT}{\bbb{T}}

\newcommand{\noind}[1]{{\setlength{\parindent}{0cm} #1}}
\newcommand{\parsk}{\par\medskip}
\newcommand{\pa}{\par\medskip}

\newcommand{\varend}{\end{document}}

\newcommand{\LP}{\left(}  \newcommand{\RP}{\right)}
\newcommand{\LQ}{\left[}  \newcommand{\RQ}{\right]}
\newcommand{\LB}{\left\{} \newcommand{\RB}{\right\}}
\newcommand{\LA}{\left\langle} \newcommand{\RA}{\right\rangle}
\newcommand{\LS}{\left|}  \newcommand{\RS}{\right|}
\newcommand{\LN}{\left\|} \newcommand{\RN}{\right\|}
\newcommand{\bLP}{\bigl(}  \newcommand{\bRP}{\bigr)}

\newcommand{\ts}{{\otimes}} \newcommand{\TS}{{\bigotimes}}

\newcommand{\es}{\emptyset}
\newcommand{\sm}{\setminus} \newcommand{\st}{\subset}
\newcommand{\eps}{\varepsilon}

\newcommand{\beware}{~$\Psi\!\Psi\!\Psi$~}

\newcommand{\EM}{\em\bf}

\newcommand{\BE}{\begin{equation}}  \newcommand{\EE}{\end{equation}}
\newcommand{\BER}{$$\begin{array}}  \newcommand{\EER}{\end{array}$$}
\newcommand{\head}[1]{
             \par\bigskip\noind{{\large\bf#1}\nopagebreak\par}\medskip}
\newcommand{\chead}[1]{\begin{center}
              \par\bigskip{\large\bf#1}\nopagebreak\par\medskip
              \end{center}}

\newcommand{\fillitem}{\L{\embox{}}\vspace{-.6cm}}

\newcommand{\dfrac}[2]{\displaystyle{\frac{#1}{#2}}}
\newcommand{\toprove}[1]{{\buildrel\mbox{?}\over#1}}
\newcommand{\DEF}{{\buildrel\mbox{\scriptsize\,\,def\,\,}\over=}}
\newcommand{\BAR}[1]{\overline{#1}}
\newcommand{\I}{\infty}
\newcommand{\al}{\alpha} \newcommand{\be}{\beta}
\newcommand{\ga}{\gamma} \newcommand{\de}{\delta}
\newcommand{\DE}{\Delta}
\newcommand{\OM}{\Omega} \newcommand{\om}{\omega}
\newcommand{\la}{\lambda}
\newcommand{\mod}{\mbox{ mod }}
\newcommand{\half}{\frac{1}{2}}
\newcommand{\NI}{{n\to\infty}}
\newcommand{\ang}{{<\!\!\!)}}
\newcommand{\arc}[1]{\breve{#1}}

\newenvironment{txeqn}[1]{\begin{equation}\label{#1}\begin{array}{l}}{
                         \end{array}\end{equation}}




\newcommand{\fn}{function}
\newcommand{\cts}{continuous}
\newcommand{\cp}{compact}
\newcommand{\tp}{topolog}
\newcommand{\gp}{group}
\newcommand{\nbhd}{neighborhood}
\newcommand{\bdd}{bounded}
\newcommand{\Bh}{Banach}
\newcommand{\sgp}{semigroup}
\newcommand{\usc}{upper semicontinuous}
\newcommand{\lsc}{lower semicontinuous}

\newcommand{\MD}{\Sigma}

\newcommand{\ob}[1]{\widetilde{#1}}




\title{Weakly Compact ``Matrices'', Fubini-Like Property and Extension of
Densely Defined Semi\gp s of Operators}
\author{Eliahu Levy\\
Department of Mathematics\\
Technion -- Israel Institute of Technology,
Haifa 32000, Israel\\
email: eliahu@techunix.technion.ac.il}


\date{}


\maketitle
\begin{abstract}
Taking {\em matrix} as a synonym for a numerical \fn\ on the Cartesian
product of two (in general, infinite) sets, a simple purely algebraic
``reciprocity property'' says that the set of rows spans a finite-dimensional
space iff the set of columns does so. Similar \tp ical reciprocity
properties serve to define strongly \cp\ and weakly \cp\ matrices,
featured in the well-known basic facts about almost periodic \fn s on \gp s
and about \cp\ operators. Some properties, especially for the weak compact
case, are investigated, such as the connection with the matrix having a
Fubini-like property for general finitely additive means. These are applied to
prove possibility of extension to the entire \sgp\ of \bdd\ densely defined \sgp s
of operators in a \Bh\ space with weak continuity properties.
\end{abstract}


\chead{Acknowledgements}

I thank Orr Shalit for introducing me to the problem of
extending densely defined \sgp s and for several discussions.
\pa

Michael Cwikel has drawn my attention to \cite{bartle} and \cite{mujica} where
similar issues are dealt with from a somewhat different point of view.
\pa

\section[Preface]{Preface: Reciprocity Between Rows and Columns, Strongly Compact
``Matrices''.}

Taking {\em matrix} as a synonym for a numerical (i.e. real or complex) \fn\
on the Cartesian product of two (in general, infinite) sets, a simple purely
algebraic ``reciprocity property'' says that {\em the set of rows spans a
finite-dimensional space iff the set of columns does so}. It might be
illuminating to sketch a proof:
If the rows span a finite-dim space, there is a finite set $F$ of them so that
all others are linear combinations of them. Otherwise put: all entries of a
column are fixed (i.e.\ the same for all columns) linear combinations of the
entries at $F$. But since $F$ is finite, the entries at $F$ of all the columns
surely span a finite-dim space, i.e.\ all depend linearly on a finite number
of them, which carries over to the entire columns by the above fixed linear
combinations. 
\pa

Assuming the matrix is bounded and viewing the rows and columns as vectors in
the $\ell^\infty$ spaces, well-known analogous \fn al-analytic reciprocity
properties obtain similarly. Thus we have: {\em the rows form a strongly
relatively \cp\ (i.e.\ precompact) set in the $\sup$-norm iff the columns
do so} (an assertion which specializes, of course, to well-known basic facts
about almost periodic \fn s on \gp s and about \cp\ operators).
To give an argument similar to the above algebraic case, suppose the rows
form a pre\cp\ set. For any $\eps>0$, there is a finite set $F$ of rows
which is an $\eps$-net (i.e.\ any row has an $\ell^\I$ distance no more
than $\eps$ from some one of them). This means that for the columns, two columns
that differ no more than some $\de>0$ at the $F$-entries differ no more than
$\de+2\eps$ at all entries. The $F$ pieces of all the columns are a bounded
set in a finite-dim space, hence precompact, i.e.\ contains a $\de$-net for
any $\de>0$ which by the above will be a $\de+2\eps$-net for the entire
columns.
\pa

Let us call a matrix where the rows, equivalently the columns, form a strongly
relatively \cp\ set in the $\sup$-norm a {\EM strongly \cp\ matrix}.
One easily obtains some properties of them (which give well-known properties
for the case of almost-periodic functions).

\begin{Prop}
The set of strongly \cp\ matrices is a closed subspace of the $\ell^\I$
space on the Cartesian product.
\end{Prop}

\begin{Prop}
\label{Prop:assc}

Let $I$, $J$ and $L$ be sets and let $f$ be a \bdd\ complex \fn\ on
$I\times J\times L$ so that

\begin{itemize}

\item $f$ is strongly \cp\ as a matrix on $L$ and $I\times J$

\item for any fixed $a\in L$, the $a$-section of $f$
$f(a,\bullet,\bullet)$ is strongly \cp\ as a matrix on $I$ and $J$

\end{itemize}

then $f$ is strongly \cp\ as a matrix on $J$ and $I\times L$ and as
a matrix on $I$ and $J\times L$.
\pa

Consequently, if a \bdd\ complex \fn\ on $I\times J\times L$ is strongly
\cp\ both as a matrix on $L$ and $I\times J$ and as a matrix on $J$
and $I\times L$, then it is so also as a matrix on $I$ and $J\times L$.
(``associativity property'').
\end{Prop}

\begin{Prop}
Suppose $I$ and $J$ are \tp ical spaces. Then if a \bdd\ complex separately
\cts\ \fn\ on $I\times J$ is strongly \cp\ as a matrix, then it is
jointly \cts. If $I$ and $J$ are compact, the converse holds: if $f$ is
jointly \cts\ it is strongly \cp\ as a matrix.
\end{Prop}

\begin{Rmk}
A \cts\ complex \fn\ $f$ on a \tp ical \gp\ $G$ is called (strongly)
almost periodic if the matrix $x,y\mapsto f(xy)$, $x,y\in G$ is strongly compact.
This clearly implies that the matrices $x,(y,z)\mapsto f(xyz)$ and
$(x,y),z\mapsto f(xyz)$ are strongly compact, hence, by Prop.\ \ref{Prop:assc}
that $(x,z),y\mapsto f(xyz)$ is strongly compact, which in turn implies that the
original matrix $x,y\mapsto f(xy)$ s strongly compact. Hence
$(x,z),y\mapsto f(xyz)$ strongly \cp\ is an equivalent definition for
almost-periodicity. \pa

Let $X,Y$ be a \Bh\ spaces and let $Y^*$ be the dual of $Y$. Let $A:X\to Y$
be a \bdd\ linear operator. Consider the matrix
$x,y^*\mapsto \LA y^*,Ax\RA$, $x$, $y^*$ in the closed unit balls of $X$ and
$Y^*$ resp. Then its rows form a relatively \cp\ set in the $\ell^\I$ space
iff $A$ is a \cp\ operator, and its columns do so iff $A^*$ is such. Thus the
matrix is strongly \cp\ $\Leftrightarrow$ $A$ is a \cp\ operator
$\Leftrightarrow$ $A^*$ is a \cp\ operator.
\end{Rmk}

\section{Weakly Compact Matrices and Fubini-Like Property} 

Now replace strong (relative) compactness by weak (relative) compactness
w.r.t.\ the $\ell^\I$ spaces. We shall need the following simple

\begin{Prop} 
\label{Prop:wkcp}
Let $I$ and $J$ be Hausdorff \tp ical spaces and let $S\st I$ be dense.
Suppose $f$ is a \bdd\ complex matrix on $I\times J$ so that all columns are
\cts\ \fn s on $I$ and all rows $f(s,\bullet)$ with {\em $s\in S$} are \cts\
\fn s on $J$. Suppose further that the set of these rows $f(s,\bullet)$,
$s\in S$ is relatively weakly \cp\ in $\ell^\I(J)$ (or, equivalently, in its
closed subspace $\cC_b(J)$ -- \bdd\ \cts\ \fn s). Then $f$ is separately
\cts\ (i.e.\ {\em all} rows $f(s,\bullet)$, $s\in I$ are \cts).
\end{Prop}

\begin{Prf}
The set of rows $f(s,\bullet)$, $s\in S$, is relatively weakly \cp\ in the
space $\cC_b(J)$. Thus, for any $r\in I$, the rows $f(s,\bullet)$ for $s\to r$
have a weak cluster point in $\cC_b(J)$ (to which they will weakly converge if
we take a finer filter, alternatively take a subnet). Anyhow, these rows tend
pointwise to $f(r,\bullet)$ since the columns are \cts. Therefore the only
possible weak cluster point is $f(r,\bullet)$, hence $f(r,\bullet)$ must be in
$\cC_b(J)$, i.e.\ \cts, and we are done.
\end{Prf}

\begin{Thm}
\label{Thm:Fubini}

Let $f$ be a \bdd\ complex matrix on the sets $I$ and $J$. TFAE:

\begin{itemize}

\item[(i)] Fubini-like property: if $M$ is any mean on $\ell^\infty(I)$
(recall this means a positive linear \fn al mapping the constant $1$ to $1$,
equivalently a finitely additive ``measure'' on $I$. $I$ is embedded into the
set of means as $\de$-``measures'') and $N$ is any mean on $\ell^\I(J)$, then
taking an ``iterated integral'' of $f$ w.r.t.\ these means does not depend on
the ``order of integration''.

\item[(ii)] $f$ can be extended to a separately \cts\ affine \fn\ 
(=matrix) by embedding $I$ and $J$ (as sets) in some $\ob{I}$ and $\ob{J}$,
\cp\ convex subsets of Hausdorff locally convex \tp ical linear spaces.
\pa

\item[(iii)] The set of rows $f(i,\bullet)$ is weakly relatively \cp\ in
$\ell^\I(J)$

\item[(iv)] The set of columns $f(\bullet,j)$ is weakly relatively compact
in $\ell^\I(I)$

\end{itemize}

In this case, we naturally say that the matrix is {\EM weakly compact}.
\end{Thm}

\begin{Prf}
(i) $\Rightarrow$ (ii):
Having the Fubini-like property, the ``iterated integral'' defines a matrix
$\ob{f}(M,N)$, which is clearly a separately \cts\ affine \fn\ on
the sets of means, these being convex \cp\ w.r.t.\ the $w^*$ \tp y from
the $\ell^\I$ spaces. $\ob{f}$ is an extension of the original $f$, when $I$
and $J$ are embedded in the sets of means as $\de$-``measures''. \pa

(ii) $\Rightarrow$ (iii) and (iv):
It clearly suffices to prove the weak relatively compactness when $I$ and $J$
themselves are convex \cp\ and the matrix $f$ separately \cts\ affine.
But note that for the $\sup$-norm space $\cA$ of \cts\ affine \fn s on a
\cp\ convex set any element of the dual is a linear combination of two
evaluation \fn als, hence the weak \tp y is the pointwise \tp y. Now each
$i\in I$ maps to the row $f(i,\bullet)$, thus mapping a \cp\ set \cts ly
w.r.t.\ to the pointwise = weak \tp y, so the set of rows is weakly compact
in the $\cA$-space, hence in the $\ell^\I$ space. Similarly for the columns.
\pa

(ii) $\Rightarrow$ (i): in (ii), one may assume that $\ob{I}$ is the closed
convex hull of the embedded $I$ and similarly for $\ob{J}$. This makes the
norm in $\ell^\I(I)$ identical with that of $\cA(\ob{I})$ (the $\sup$-norm
space of \cts\ affine \fn s), so $\cA(\ob{I})$
may be viewed as a subspace of $\ell^\I(I)$, similarly for $J$. Any mean
on $\ell^\I(I)$ ``collapses'' on $\cA(\ob{I})$ to an evaluation \fn al at a
point of $\ob{I}$, similarly for $\ob{J}$. For evaluation \fn als Fubini
is immediate. \pa
 
(iii) $\Rightarrow$ (ii): We wish to invoke Prop.\ \ref{Prop:wkcp}. Indeed,
taking ``iterated integral'' first on the first argument then on the second,
one extends $f$ to $\ob{f}(M,N)$ on the \cp\ convex sets of means, \cts\
(and affine) in $N$ (i.e.\ with \cts\ rows), and also \cts\ (and affine) in
$M$ if $N=\de_j$, $j\in J$ hence if $N$ is a finite convex combination of
$\de_j$'s, these being dense in the set of means on $J$. Also the sup-norm
for \cts\ affine \fn s on the set of means on $J$, such as the rows of
$\ob{f}$, coincides with that of their restriction to $J$, making $\ell^\I(J)$
isometric to the sup-norm space of these \cts\ affine \fn s, and our assumption
that the rows form a weakly relatively \cp\ set translates accordingly.
Therefore by Prop.\ \ref{Prop:wkcp} $\ob{f}$ is separately \cts\ (and affine)
and we have (ii).
\end{Prf}

Thus the matrix being weakly \cp\ is determined by equations -- the
equality of the ``iterated integrals'', and one concludes

\begin{Cor}
The set of weakly \cp\ matrices is a closed subspace of the $\ell^\I$
space on the Cartesian product. 
\end{Cor}

While jointly \cts\ matrices on two \cp\ spaces are strongly compact,
separately \cts\ matrices on two \cp\ spaces are weakly compact:
(Thm.\ \ref{Thm:sepcomp} is not needed in \S\ref{sec:dens})

\begin{Thm}
\label{Thm:sepcomp}

Let $I$ and $J$ be \cp\ Hausdorff spaces, and let $f$ be a \bdd\ complex
matrix on $I\times J$, separately \cts. Then the matrix $f$ is weakly compact.
\end{Thm}

\begin{Prf}
{}

\begin{itemize}

\item Assume $J$ is convex \cp\ and every $f(i,\bullet)$ is
affine. Then for the $\sup$-norm space $\cA(J)$ of \cts\ affine \fn s on
$J$ every element of the dual is a linear combination of two evaluation
\fn als, hence the weak \tp y is the pointwise \tp y.
$i\mapsto f(i,\bullet)$ maps the \cp\ $I$ \cts ly to $\cA(J)$ with the
weak = pointwise \tp y, hence the image is weakly \cp\ so the rows
form a weakly \cp\ set.

\item Assume $I$ metrizable. Let $\Pr(J)$ be the \cp\ convex
space of the probability measures on $J$ (with the $w^*$-\tp y
from the \cts\ \fn s). Extend $f$ to a matrix $\ob{f}$ on $I\times\Pr(J)$ by
$\ob{f}(i,\mu):=\int f(i,j)\,d\mu(j)$. $\ob{f}$ is \cts\ in $\mu$, but also
in $i$ -- on the metrizable $I$ continuity of \fn s is determined by
sequences and one uses Lebesgue's convergence theorem. So by the previous
$\bullet$, $\ob{f}$ is weakly compact, hence so is $f$.

\item $I$ and $J$ general. It is a well-known fact (due to Eberlein and
\v{S}mulian) that a subset $A$ of a \Bh\ space, such as
$\ell^\I$, is weakly relatively \cp\ iff every countable subset of $A$ is
so. Therefore to prove that the set of rows $f(i,\bullet)$ is weakly
relatively \cp\ it suffices to prove that every sequence of rows
$f(i_n,\bullet)$ is such. But such a sequence may be viewed as a \cts\
mapping from $J$ with metrizable \cp\ image $Q\st\bC^{\bN}$. \pa

Consider $\be\bN$ -- the Stone-\v{C}ech compactification of $\bN$, i.e.\ the
\cp\ space of the ultrafilters in $\bN$. One can define a matrix $m$ on
$\be\bN\times Q$ by $m(\tau,q):=\lim_{n\to\tau}q_n$. $m$ is \cts\ in $\tau$,
but also in $q$, because substituting $q(j)=\LP f(i_n,j)\RP_{n=1,2,\ldots}$
gives

$$m(\tau,q)=\lim_{n\to\tau}f(i_n,j)=f\LP\lim_{n\to\tau}i_n,j\RP,$$

where one uses the continuity of f in i and the compactness of $I$, ensuring
the existence of $\lim_{n\to\tau}i_n$. \pa

So by the previous $\bullet$, $m$ is weakly compact, hence so is the sub-matrix of $f$
formed by the sequence of rows, which is a substitution in $m$.

\end{itemize}
\end{Prf} 

\begin{Thm} 
\label{Thm:cts}

Let $A$ be \cp\ convex (subset of a Hausdorff locally convex \tp ical
linear space), let $I$ be a Hausdorff \tp ical space and let $S\subset I$ be
dense. Suppose $f(\al,s)$ is a complex \fn\ on $A\times I$ satisfying: \pa

\begin{itemize}

\item[(i)] $f$ is \cts\ affine in $\al$ on $A$ for each fixed $s\in S$.

\item[(ii)] $f$ is \cts\ in $s$ on $I$ for each fixed $\al\in A$.

\end{itemize}

Then $f$ is \cts\ in $\al$ on $A$ also for each fixed $r\in I$, i.e.\ it
is separately \cts\ on $A\times I$.
\end{Thm}

\begin{Prf}
By Prop.\ \ref{Prop:wkcp} it is enough to prove that the matrix $f$ on
$A\times S$ is weakly compact. By Thm.\ \ref{Thm:Fubini} one just needs to
establish the Fubini-like property for $f$ w.r.t.\ any means on $\ell^\I(A)$
and $\ell^\I(S)$.
But any mean on $\ell^\I(A)$ ``collapses'', on the sections $f(\bullet,s)$
$s\in S$ parallel to $A$ (the columns), to an evaluation at some point
$\al\in A$.
Thus, if $M$ is a mean on $\ell^\I(S)$, then both iterated integrals give
$M(f(\al,\bullet))$. So Fubini holds and we are done.
\end{Prf}

\begin{Rmk}
The analogue of the ``associativity property'' Prop.\ \ref{Prop:assc} does
not hold for weakly \cp\ matrices. An example is $B$ = the unit ball
of a Hilbert space $H$, $L$ = the set of operators in $H$ of norm $\le1$
and $f$ defined on $B\times B\times L$ by $f(x,y,A):=\LA Ax,y\RA$. It is
weakly \cp\ as a matrix on $(x,A)$ and $y$ and as a matrix on $(y,A)$ and
$x$ but not as a matrix on $(x,y)$ and $A$.
\end{Rmk}

\section{A Digression: Upper- and Lower- Semi\cts\ Envelopes}
\label{sec:sc}

Let $S$ be a dense subset of a completely regular \tp ical space $I$. For
any \bdd\ real \fn\ $f$ defined on $S$ let $f^-$ be defined (on the whole $I$)
as the pointwise infimum of all \cts\ \fn s (on $I$) that majorize $f$ on the
$S$. $f^-$ is \usc\ (an infimum of \cts, even of \usc\ \fn s is always \usc,
and any \usc\ \fn\ is the infimum of all \cts\ \fn s which majorize it) and is
the smallest \usc\ \fn\ majorizing $f$ on $S$. \pa

Note that $f^-$ at a point $r\in I$ is, in fact, defined locally: is the same
for two \fn s $f$ coinciding on $S$ in a \nbhd\ of $r$. \pa

Equivalently, since for any function $f$ on $S$ the function of $t$ in $I$
$\limsup_{s \to t}f(s)$ is \usc, we have $f^-(t)=\limsup_{s\to t}f(s)$. So this
is another way to define $f^-$. Similarly for $f_-$ below. \pa

Replacing majorizing by minorizing (thus upper by \lsc) one gets $f_-$. Clearly
$f^-\ge f_-$ on $I$. \pa

Surely, $f$ can be extended to a \cts\ \fn\ on $I$ iff $f^- = f_-$ everywhere
on $I$. \pa

\section{Densely Defined Semi\gp s of Operators}
\label{sec:dens}

\begin{Thm}
Let $\MD$ be a \tp ical \sgp\ which is an open sub-\sgp\ of a \tp ical \gp\
which is a Baire space (e.g.\ a locally \cp\ \gp\ or a Polish \gp).
Let $\cF$ be a filter in $\MD$ having a basis consisting of open sets and
converging to the unit element of the \gp\
(think of $\MD=(0,\I)$ with addition and $\cF$ -- tending to $0^+$).
Let $S$ be a dense sub-\sgp\ of $\MD$. Let $X$ be a \Bh\
space, taken as real. (Elements of $X$ denote by $x$ and elements of its dual
$X^*$ by $\rho$). Let $T_s, s\in S$ be a uniformly \bdd\ \sgp\ of linear
operators on $X$ defined on $S$. Suppose $T_s\to 1$ weakly as $s$ tends
to $\cF$ in $S$. Then $T_s, s\in S$ can be extended to a weakly \cts\ (in $r$)
\sgp\ $T_r, r\in\MD$.
\end{Thm}

\begin{prf}{\em
Replacing the norm in $X$ by the equivalent norm
$$\|x\|_1:=\sup_s\|T_s(x)\|$$
we may and do assume that the $T_s$ are contractions. \pa

Let $B^*$ be the unit ball of the dual $X^*$, with the $w^*$-\tp y.
The members $x\in X$ may be identified with the \fn s on $B^*$\quad
$x(\rho):=\rho(x)$,\,\,$\rho\in B^*$.
(And denote that \fn\ by the same $x$.)
This is an isometry of $X$ into $\ell^\infty(B^*)$.
Moreover, the image of $X$ are precisely all \fn s on $B^*$ that are \cts\
affine and vanish at $0$, i.e.\ extend to linear on $X^*$.
(That follows from a well-known fact about \Bh\ spaces, due to Dieudonn\'e:
any linear \fn al on $X^*$ which is \cts\ on $B^*$ w.r.t.\ the $w^*$ \tp y
comes from an element of $X$.) \pa

The action of $T_s$ on $X$ translates as follows:
$$(T_sx)(\rho)=\rho(T_sx)=(T_s^*\rho)(x)=x(T_s^*\rho).$$
Thus our \sgp\ translates to a \sgp\ induced by a \sgp\ of 
(\cts\ affine fixing $0$) transformations $T_s^*$ of the space $B^*$. \pa

Our task, extending the \sgp, will be accomplished if, for each $x$, we can
extend the \fn\ $(\rho,s)\mapsto x(T_s^*\rho)$ to a \fn\ on
$(\rho,t)\in B^*\times\MD$ which is {\em separately \cts}. Indeed, then by
continuity it will be affine in $\rho$, vanishing at $0$, hence by
Dieudonn\'e's theorem will be of the form $(T_tx)(\rho)$, thus defining
$T_tx$ and $T_t$, $t\in\MD$, which will depend weakly \cts ly on $t$ and will
coincide with the original $T_s$ for $s\in S$, therefore, by continuity, will
satisfy the \sgp\ identity. \pa

But note, that if we can extend each $s\mapsto x(T_s^*\rho)$, for fixed
$x$ and $\rho$, to a \fn\ \cts\ on $\MD$, then Thm.\ \ref{Thm:cts} (taking
there $A=B^*$ and $I=\MD$) will insure the separate continuity. \pa

Thus we are left with proving that $s\mapsto x(T_s^*\rho)$, for fixed
$x$ and $\rho$, can be \cts ly extended to $\MD$. \pa

For each $\rho\in B^*$ and $h\in\ell^\infty(B^*)$ there is the ``orbit \fn''
on $S$\quad$O_\rho(h):=s\mapsto h(T_s^*(\rho))$. The map $h\mapsto O_\rho(h)$
is, of course, a \bdd\ linear operator from $\ell^\infty(B^*)$ to
$\ell^\infty(S)$. Moreover, the action $T_sh(\rho):=h(T_s^*\rho)$ translates
to
\begin{eqnarray*}
&&O_\rho(T_sh)(s')=T_sh(T_{s'}^*(\rho))=h(T_s^*T_{s'}^*\rho)=\\
&&=h((T_{s'}T_s)^*\rho)=h(T_{s's}^*\rho)=(O_\rho(h))(s's),
\end{eqnarray*}
that is, to the shift $f\mapsto(s'\mapsto f(s's))=f(\bullet\cdot s)$. \pa

Now, the fact that $T_s\to1$ as $s\to\cF$ in $S$ translates, for the
orbit \fn\ $x(T_s^*\rho)$, which is the \fn\ that we have to \cts ly
extend to $\MD$, to ``$f(\bullet\cdot s)$ tends to $f$ weakly in $\ell^\I(S)$
as $s\to\cF$ in $S$''. \pa

Hence it suffices to prove:
}\end{prf}

\begin{Lem}
Let $f$ be a real \bdd\ \fn\ on $S$ with the property: \pa
 
$s'\mapsto f(s's)$, i.e.\ $f(\bullet\cdot s)$, tends to $f$ weakly
in $\ell^\I(S)$ as $s\to\cF$ in $S$. \pa

Then $f$ can be extended to a \cts\ \fn\ on $\MD$.
\end{Lem}

\begin{Prf}
We shall use $f^-$ and $f_-$ defined in \S\ref{sec:sc}. Our task is to prove
that $f^-=f_-$ everywhere on $\MD$. \pa
 
As $s\to\cF$ in $S$, $f(\bullet\cdot s)$ tends to $f$ weakly in $\ell^\I(S)$,
hence, in particular, $f(s's)\to f(s')$ pointwise (i.e.\ for each $s'\in S$).
But this implies that $f^-(s's)\to f(s')$ and $f_-(s's)\to f(s')$ when
$s\in\MD$ tends to $\cF$ in $\MD$. Indeed, fix $s'$ and $\eps>0$. Choosing
$U$ small enough in an open basis of $\cF$ one has
$f(s')-\eps\le f(s's)\le f(s')+\eps$
for $s\in U\cap S$, so the \usc\ \fn\ $g$, defined as the constant $f(s')+\eps$
in the open $s'U$ and as $+\infty$ elsewhere, majorizes $f$ on $S$, hence
majorizes $f^-$ everywhere, making $f^-(s's)\le f(s')+\eps$ for $s\in U$.
In a similar manner one bounds $f_-$ below using a \lsc\ \fn\ equal to
$f(t)-\eps$ in $s'U$ and to $-\infty$ elsewhere. \pa

Consider {\em means} $M$ on $\ell^\infty(S)$ (i.e.\ positive linear
\fn als mapping the constant $1$ to $1$) such that for restrictions of \fn s
\cts\ on $\MD$, $M$ evaluates the \fn\ at some fixed $r\in\MD$. They are
characterized by
\BE\label{eq:Ma} 
g^-(r)\ge M(g)\ge g_-(r)\quad\mbox{for }g\in\ell^\I(S).
\EE
(A particular case is evaluation \fn al at some $r\in S$.) \pa

We have $f(\bullet\cdot s)\to f$ weakly in $\ell^\I(S)$ as $s\to\cF$ in $S$,
therefore
\BE\label{eq:Mb}
M(f(\bullet\cdot s))\to M(f)\quad\mbox{as }s\to\cF\mbox{ in }S
\EE

Plug in (\ref{eq:Ma}) $g=f(\bullet\cdot s)$. Note that then 
$g^-=f^-(\bullet\cdot s)$, $g_-=f_-(\bullet\cdot s)$
(this follows from $\MD$ being an open part of a \tp ical \gp\ where shifts
are, of course, homeomorphisms). One gets
\BE\label{eq:Mc}
f^-(rs)\ge M(f(\bullet\cdot s))\ge f_-(rs),
\EE
and (\ref{eq:Mb}) and (\ref{eq:Mc}) imply:
\BE\label{eq:Md}
\liminf_{s\to\cF\mbox{ in }S}f^-(rs)\ge
M(f)\ge\limsup_{s\to\cF\mbox{ in }S}f_-(rs).
\EE
Now, for any $f^-(r)\ge a\ge f_-(r)$, using Hahn-\Bh\ one gets a mean as
above so that $M(f)=a$.
(Just extend to a positive linear \fn al the linear \fn al on the space
generated by $f$ and the restrictions of \fn s \bdd\ and \cts\ on $\MD$,
which evaluates the latter at $r$ and gives value $a$ for $f$.)
Thus (\ref{eq:Md}) reads:
$$\liminf_{s\to\cF\mbox{ in }S}f^-(rs)\ge f^-(r),\qquad
f_-(r)\ge\limsup_{s\to\cF\mbox{ in }S}f_-(rs).$$
By the semicontinuity of $f^-$ and $f_-$ one can replace here
$s\to\cF\mbox{ in }S$ by $s\to\cF\mbox{ in }\MD$.
Indeed, if $\la<f^-(r)$ we know that the set $\{s\,|\,f^-(rs)\ge\la\}$ is
closed and contains $S\cap U$ for some open $U$ in $\cF$, therefore contains
$U$. Similarly for $f_-$.
Also, since $f^-$ is \usc,
$$f^-(r)\ge\limsup_{s\to\cF\mbox{ in }\MD}f^-(rs),$$
making $f^-(rs)\to f^-(r)$ and similarly $f_-(rs)\to f_-(r)$ as $s\to\cF$ in
$\MD$. \pa

Since, by what we found, for $r=s'\in S$ these limits are $f(s')$, one
concludes that $f^-$ and $f_-$ coincide with $f$ on $S$. \pa

Now, $f^- - f_-$ is \usc\ and its value at $s's$ tends to $0$ as $s\to\cF$ in
$\MD$ for each $s'\in S$. Hence for any $\eps>0$, the open set
$\{r\in\MD\,|\,f^-(r)-f_-(r)<\eps\}$ is dense in $\MD$. Thus by Baire category
the set $\{r\in\MD\,|\,f^-(r)=f_-(r)\}$ is dense $G_\delta$.
With $f^-(rs)\to f^-(r)$ and $f_-(rs)\to f_-(r)$ as $s\to\cF$ in $\MD$, this
makes $f^-=f_-$ everywhere, which concludes the proof of the lemma and the
theorem.
\end{Prf}

\end{document}